\documentclass[]{article}
 
\usepackage[modulo,switch]{lineno}
\modulolinenumbers[1]
\usepackage[utf8]{inputenc}
\usepackage{amsmath}
\usepackage{amssymb}
\usepackage{graphicx}
\usepackage{xcolor}
\usepackage{dashrule}
\usepackage{subfig}
\usepackage[above,below]{placeins}

\makeatletter\@ifundefined{date}{}{\date{}}
\makeatother

\evensidemargin\oddsidemargin \hoffset-29.4mm \voffset-28.4mm
\usepackage{bm}
\usepackage{hyperref}

\newcommand{\etal}{\textit{et~al.}}
\renewcommand{\vec}[1]{\bm{#1}}

\newcommand{\dd}{\mathrm{d}}
\newcommand{\pd}{\partial}

\newcommand{\pfrac}[2]{\frac{\pd #1}{\pd #2}}

\newcommand{\grad}{\nabla}
\renewcommand{\div}{\nabla\cdot}

\renewcommand{\exp}[1]{\mathrm{e}^{#1}}

\newcommand{\unit}[1]{\ \mathrm{#1}}

\usepackage[only,llbracket,rrbracket]{stmaryrd}
\newcommand{\DGaverage}[1]{\left\{ \!\! \left\{ #1 \right\} \!\! \right\}}
\newcommand{\DGjjump}[1]{\left\llbracket #1 \right\rrbracket}
\newcommand{\DGjump}[1]{\left[#1\right]}

\usepackage[normalem]{ulem}

\usepackage{nicefrac}

\renewcommand{\varepsilon}{\mathcal{E}}
\usepackage{booktabs}

\usepackage[singlespacing]{setspace}
\begin{document}

\title{High-Order Non-Conforming Discontinuous Galerkin Methods for the Acoustic Conservation Equations}
\author{Johannes Heinz\footnote{Institute of Mechanics and Mechatronics, TU Wien, Vienna, Austria (\texttt{johannes.heinz@tuwien.ac.at})}, Peter Munch\footnote{High-Performance Scientific Computing, University of Augsburg, Augsburg, Germany (\texttt{peter.muench@uni-a.de})}, Manfred Kaltenbacher\footnote{Institute of Fundamentals and Theory in Electrical Engineering (IGTE), TU Graz, Graz, Austria (\texttt{manfred.kaltenbacher@tugraz.at})}}

\maketitle\thispagestyle{empty}
\begin{abstract}
This work compares two Nitsche-type approaches to treat non-conforming triangulations for a high-order discontinuous Galerkin (DG) solver for the acoustic conservation equations.
The first approach (point-to-point interpolation) uses inexact integration with quadrature points prescribed by a primary element.
The second approach uses exact integration (mortaring) by choosing quadratures depending on the intersection between non-conforming elements.
In literature, some excellent properties regarding performance and ease of implementation are reported for point-to-point interpolation.
However, we show that this approach can not safely be used for DG discretizations of the acoustic conservation equations since,
in our setting, it yields spurious oscillations that lead to instabilities.
This work presents a test case in that we can observe the instabilities and shows that exact integration is required to maintain a stable method.
Additionally, we provide a detailed analysis of the method with exact integration.
We show optimal spatial convergence rates globally and in each mesh region separately.
The method is constructed such that it can natively treat overlaps between elements.
Finally, we highlight the benefits of non-conforming discretizations in acoustic computations by a numerical test case with different fluids.
\end{abstract}

\section{Introduction}
The main benefit of non-conforming interfaces (NCIs) is the ability to handle arbitrary element connections.
In acoustic simulations, we require different element sizes in different regions of a triangulation, e.g., due to wave propagation through inhomogeneous media.
NCIs can realize the jump in element sizes without the use of transition regions which usually contain strongly distorted elements~\cite{Quiroz1995,Kaltenbacher2015}.
This way, it is possible to reduce degrees of freedom (DoFs) needed without introducing errors related to elements with bad quality~\cite{Flemisch2006,Kaltenbacher2015}.
Additionally, algorithms that can handle NCIs can simplify mesh generation since it is possible to generate the grids in a modular way~\cite{Quiroz1995}.

Overlapping elements further reduce the difficulties in mesh generation since they can be constructed without paying attention to adjacent regions at all.
One famous example of this is the overset grid method~\cite{Benek1985};
a structural mesh serves as the background mesh in which complex geometries can be embedded.
This is done by overlaying the corresponding meshes and deleting the elements of the background mesh that completely overlap the embedded mesh.

Besides mentioned advantages of NCIs, some applications, like a rotating fan, require NCIs.
To compute the aeroacoustic sound field, we need two mesh regions, one of which is rotating, see e.g.~\cite{Schoder2020}.
A conforming mesh at each time step is not feasible; 
using NCIs in such applications is the obvious solution.
However, this requires the non-conforming interface to lie precisely on top of each other, which is only possible using curved elements, cf.~\cite{Pezzano2022}.
A slightly different approach is to use methods that can also handle element overlaps between the triangulations.
This way, the fixed and rotating domain can still pick values for the fluxes at the overlapping boundaries, with the difference that these values are defined inside elements of the other triangulation.

There exist three different ways to handle non-conformities.
The most common method in literature is the Mortar method, first introduced by Bernadi~\etal~\cite{Bernardi1990}.
The Mortar method is a projection-based method that typically uses Lagrange multipliers to enforce coupling; this requires additional DoFs at the interface.
Coupling of the second order wave equation using Mortar methods has been successfully applied in~\cite{Flemisch2006}.

Another way to couple non-conforming meshes is through interpolation-based methods, such as INTERNODES (INTERpolation for NOn-conforming DEcompositionS)~\cite{Deparis2016}.

Nitsche~\cite{Nitsche1971} presented the idea of including Dirichlet boundary conditions (DBCs) in the weak form. 
Methods using this idea are consequently named Nitsche-type methods.
Discontinuous Galerkin (DG) schemes use this idea at all element boundaries already~\cite{Arnold2002}.
Therefore, we believe that using Nitsche-type methods to couple meshes via NCIs is the most natural way to couple DG schemes.
Here we can distinguish between schemes that use exact and inexact integration.
For methods with inexact integration, we use integration points dictated by elements on the NCI and evaluate needed quantities in the non-conforming attached elements.
Hermann~\etal~\cite{Hermann2010} used this approach in a two dimensional DG setting for seismic waves on meshes with possibly different element types.
Laughton~\etal~\cite{Laughton2021} refer to this method as point-to-point interpolation method.
Methods using exact integration collect integration rules on the intersections between the connected elements.
This procedure is commonly referred to as ``mortaring'', and the intersections are often called ``mortars'' (without any relation to the Mortar method).
Nitsche-type mortaring has been successfully applied in a FEM setting by e.g.~\cite{Kaltenbacher2018a, Roppert2020} for the inhomogeneous wave equation and Maxwell equations, respectively.
The procedure of mortaring is the same for Mortar and Nitsche-type mortaring methods.
The difference in both methods is that Nitsche-type methods enforce the coupling point-wise (via numerical fluxes).
On the other hand, Mortar methods enforce the coupling via an integral (using Lagrange multipliers).

Laughton~\etal~\cite{Laughton2021} compared the Nitsche-type mortaring method to the point-to-point interpolation method regarding performance and accuracy in a DG setting.
The advantage of point-to-point interpolation over methods with mortaring is its ease of implementation~\cite{Laughton2021}.
For the compressible Euler equations in two dimensions, it is shown that point-to-point interpolation outperforms the method with mortaring, considering polynomial degrees between~$3$ and~$7$~\cite{Laughton2021}.
We expect the performance to close up for long run-times on static triangulations (the quadrature rules of the intersections and the mapping of obtained integration points have to be setup only once).
However, we suspect the performance to diverge even more on moving meshes, where the intersections and the mappings have to be updated every time step or Runge--Kutta stage.
The disadvantage of point-to-point interpolation is that it introduces numerical errors related to aliasing.
Methods using mortaring do not face this issue.
To obtain similar errors for point-to-point interpolation compared to the Nitsche-type mortaring,~\cite{Laughton2021} increases the number of quadrature points.

Solving the scalar acoustic wave equation utilizing conforming finite element methods (FEM) has some unattractive peculiarities.
It requires specific time-stepping schemes and suffers from numerical dispersion (see e.g.~\cite{Cohen2002, Kaltenbacher2015}).
Transforming the acoustic wave equation to a first-order system yields the acoustic conservation equations.
These acoustic conservation equations do not include a second-order temporal derivative; 
thus, standard time-stepping methods, such as Runge--Kutta methods, can be applied.
Furthermore, the velocities of non-harmonically vibrating surfaces natively appear in governing equations, making a straightforward application of these velocities as boundary conditions (BCs) possible.
Additionally, conservation laws are ideally suited for finite volume or DG discretizations~\cite{Hesthaven2008}, and it is possible to find less dispersive schemes by adding numerical diffusion using numerical fluxes. 

We applied the point-to-point interpolation method in~\cite{Heinz2022} and showed that it provides optimal rates of convergence in space.
Later, we observed instabilities for some element configurations using this method.
Within this work, we show that for DG discretizations of the acoustic conservation equations, it is not safe to use point-to-point interpolation since the method is not only less accurate but yields spurious oscillations that lead to instabilities in some cases.
To the best of the authors' knowledge, no Nitsche-type mortaring formulation exists for the acoustic conservation laws in literature.
We present a test case in which mentioned instabilities occur and show that using exact integration via mortaring is a remedy.

Additionally, we provide in-depth convergence studies for the Nitsche-type mortaring approach and show optimal spatial convergence rates on the global computational domain and separately on domains with coarse and high resolution.

\section{Governing Equations\label{sec:goveq}}
The wave equation reads
\begin{align}
\frac{1}{c^2}\pfrac{^2 p}{t^2}-\rho\div\left(\frac{1}{\rho}\grad p\right)=f\quad&\text{in }\Omega\times[0,T],
\end{align}
on a domain~$\textstyle \Omega\subset\mathbb{R}^d$ of dimension~$\textstyle d$ and in a time interval~$\textstyle [0,T]$.
Here,~$p$ is the acoustic pressure,~$c$ is the sound speed, and the underlying material's density is $\textstyle{\rho}$.
The wave equation is a reformulation of the acoustic conservation equations of momentum and mass
\begin{align}
  \rho\pfrac{\vec{u}}{t}+\grad p=0\quad&\text{in }\Omega\times[0,T],  \label{eq:gov1}	\\
  \frac{1}{c^2}\pfrac{p}{t}+\rho\div \vec{u}=F\quad&\text{in }\Omega\times[0,T],  \label{eq:gov2}\\
  p=g_p\quad&\text{on }\pd\Omega_p^\mathrm{D},\label{eq:pdbc}\\
  \vec{u}=\vec{g}_u\quad&\text{on }\pd\Omega_u^\mathrm{D},\label{eq:udbc}\\
  \rho c\vec{u}\cdot\vec{n}=Yp\quad&\text{on }\pd\Omega^{\mathrm{Y}}.\label{eq:ybc}
\end{align}
At boundaries we can apply pressure DBCs~\eqref{eq:pdbc}, velocity DBCs~\eqref{eq:udbc}, and admittance BCs~\eqref{eq:ybc} by setting the normal component of the velocity and a certain admittance $Y$.

\section{Numerical Method}
\subsection{Notation} 

The physical domain $\Omega$ is represented by the computational domain~$\Omega_h(t)=\bigcup_{i=1}^{N_{\textrm{el}}}\Omega_{e_i}\in \mathbb{R}^d$, with the space dimension $d$.
Within this work, it consists of~$N_{\textrm{el}}$ possibly overlapping rectangular/hexahedral finite elements and is bounded by~$\Gamma_h=\pd \Omega_h$.
A finite element spans~$\Omega_{e}$ and is bounded by $\pd\Omega_{e}$.
The solution is continuous inside elements and, due to the nature of DG, discontinuous between elements.
The acoustic particle velocity $\vec{u}$ and acoustic pressure $p$ are subject to the broken polynomial spaces $\mathcal{V}_h$ for the corresponding test and trial functions
\begin{align}
    \mathcal{V}_h^u&=\Big\{\vec{u}_h\in   \ [L_2(\Omega_h)]^d:\vec{u}_h(\vec{x})|_{\Omega^{e}}=\tilde{\vec{u}}(\vec{\xi})|_{\tilde{\Omega}^{e}}\in[\mathcal{P}_{k_u}(\tilde{\Omega}^{e})]^d,\forall e\in[1,N_{\mathrm{el}}]\Big\},\\
  \mathcal{V}_h^p&=\Big\{p_h\in\ L_2(\Omega_h) :p_h|_{\Omega^{e}}=\tilde{p}_h|_{\tilde{\Omega}^{e}}\in \mathcal{P}_{k_p}(\tilde{\Omega}^{e}),\forall e\in[1,N_{\mathrm{el}}]\Big\}.
\end{align}
Here~$\textstyle \mathcal{P}_k$ is the space of polynomial functions with order~$\textstyle k$ on a reference element.
Coordinates in the physical space are~$\vec{x}=(x_1,...,x_d)^T$; their representation on a reference element are~$\vec{\xi}=(\xi_1,...,\xi_d)^T$.
To transfer between~$\vec{x}$ and~$\vec{\xi}$ a bidirectional mapping
\begin{align}
\vec{\varphi} :\left\{\begin{array}{l}
                       \Omega^e\rightarrow \tilde{\Omega}^e\\
                       \vec{x}\mapsto \vec{\xi}=\vec{\varphi}(\vec{x},\Omega^e)
                        \end{array}\right.,\quad
  \vec{\varphi}^{-1} :\left\{\begin{array}{l}
                      \tilde{\Omega}^e\rightarrow \Omega^e\\
                      \vec{\xi}\mapsto \vec{x}=\vec{\varphi}^{-1}(\vec{\xi},\Omega^e)
                        \end{array}\right.,
\end{align}
can be used.
The discrete representations of the continuous pressure and velocity fields in the reference space read
\begin{align}
	\tilde{\vec{u}}_h^e (\vec{\xi})
	=
	\sum_{i=1}^{n_{N^{k_{\vec{u}}}}} N_{i}^{k_{u}}(\vec{\xi})\vec{u}_{i}^{e},
\quad
	\tilde{p}_h^e (\vec{\xi}) 
	=
	\sum_{i=1}^{n_{N^{k_{p}}}} N_{i}^{k_p}(\vec{\xi})p_{i}^{e},
\end{align}
with the number of shape functions~$n_{N^k}$ defined on a volume element, e.g. for the pressure in the one dimensional case~$n_{N^{k_p}}=k_p+1$.
The shape functions~$N_{i}^{k}$ are constructed by Lagrange polynomials of degree~$k$.
Within this work, the same polynomial orders $k$ for velocity and pressure~($k_{\vec{u}} =k_p$) are utilized. 
We denote the interior information of an element~$\Omega^e$ with~$(\cdot)^-$ and the exterior information of adjacent elements with~$(\cdot)^+$. 
Consequently, the current element (from now on called ``primary element'') is denoted as~$\Omega^e_-$, and the neighboring elements (or ``secondary elements'') are described as~$\Omega^e_+$. 
The outward pointing normal vectors of the primary element are~$\vec{n}^-$, since facets of primary and secondary elements coincide,~$\vec{n}=\vec{n}^-=-\vec{n}^+$. 
Accordingly, any scalar or vectorial quantity~$b$ is implicitly defined on the primary element if no superscript explicitly assigns it to the primary or secondary element $b=b^-$.
We choose the notation for the averaging operator~$\DGaverage{\cdot}$, jump operator~$\DGjump{\cdot}$, and normal jump operator $\DGjjump{\cdot}$ according to Bassi~\textit{et~al.}~\cite{Bassi2006,Bassi2007}.
They are $\DGaverage{b} = \nicefrac{(b^- + b^+)}{2}$, $\DGjump{b} = b^- - b^+$, and~$\DGjjump{b} = b^- \otimes n^- + b^+ \otimes n^+$.
Hereinafter, the integrals are written in the compressed notation~$\left (a, b\right )_{\Omega_e} =\int_{\Omega_e} a \cdot b \ \dd\Omega$ and~$\left (a, b\right )_{\pd\Omega_e} = \int_{\pd\Omega_e} a\cdot b \ \dd\Gamma$, where the operator~$\cdot$ indicates an inner product and $a$ represents an arbitrary quantity of the same dimension as~$b$.
All operators are given in the notation considering element boundaries; 
therefore, each facet becomes a primary and secondary element.
For numerical integration, we employ Gaussian quadrature.
On an element of spatial dimension $d$ we use~$n_q=(k +1)^d$ quadrature points for the volume integrals and~$n_q=(k +1)^{d-1}$ on element faces.
Boundary conditions are applied using a mirror principle, cf.~\cite{Hesthaven2008}.
While pressure and velocity DBCs are defined as
\begin{alignat}{2}
  p^+=-p^-+2g_p  ;\quad \vec{u}^+&=\vec{u}^-\quad&&\text{on }\pd\Omega_p^\mathrm{D},\label{eq:pdbcval}\\
  \vec{u}^+=-\vec{u}^-+2\vec{g}_u;\quad   p^+&=p^-\quad&&\text{on }\pd\Omega_u^\mathrm{D}.\label{eq:udbcval}
\end{alignat}
Admittance BCs read
\begin{align}
\vec{u}^+ =\left(\frac{2Y} {\rho c}p^- -\vec{u}^- \cdot\vec{n}\right) \vec{n};\quad p^+= p^-\quad \text{on }\pd\Omega^{\mathrm{Y}}.\label{eq:admitt}
\end{align} 
Reflecting BCs and first-order absorbing BCs~\cite{Engquist1977} (ABC) are achieved by setting the admittance to $Y=0$ and $Y=1$, respectively.
If the first order ABC is insufficient, a corresponding perfectly matched layer (PML) formulation is provided in~\cite{Hueppe2012} for conforming FEM formulations.

\subsection{Spatial Discretization}

The numerical method, without non-conformities, has been described briefly in~\cite{Heinz2021}.
Within this section, we will recall it in a more detailed manner to be able to extend the formulation.

The semi-discrete system is obtained as usual (cf.~\cite{Hesthaven2008}).
The governing equations are multiplied by the test functions~$\textstyle \vec{w}_h$ and~$\textstyle q_h$, and integrated over the computational domain~$\textstyle \Omega_h$.
For DG schemes, it is crucial to perform the integration by parts to ensure boundary terms exist.
With this, we end up with a corresponding weak formulation.
For given equations, it is also possible to perform a second integration by parts to obtain the strong formulation, cf.~\cite{Hesthaven2008}, which is used in~\cite{Hochbruck2014}.
Eventually, numerical fluxes (denoted by the superscript~$\textstyle ^*$) are introduced into the boundary integrals.
This results in the semi-discrete system of equations
\begin{alignat}{2}
  \Bigg(\vec{w}_h,\pfrac{\vec{u}_h}{t}\Bigg)_{\Omega^e}-\left(\frac{1}{\rho}\div\vec{w}_h, p_h\right)_{\Omega^e} +\left(\frac{1}{\rho}\vec{w}_h\cdot\vec{n}, p_h^*\right)_{\pd\Omega^e}&=0\quad &&\forall\vec{w}_h\in\mathcal{V}_h^u,\\
  \bigg(q_h,\pfrac{p_h}{t}\bigg)_{\Omega^e}-\left(\rho c^2\grad q_h,\vec{u}_h\right)_{\Omega^e}+\left(\rho c^2q_h\cdot\vec{n}, \vec{u}_h^*\right)_{\pd\Omega^e}&=\left(c^2q_h,f\right)_{\Omega^e}\quad &&\forall q_h\in\mathcal{V}_h^p.	
\end{alignat}
We use Lax--Friedrichs fluxes, as also done in \cite{Nguyen2011, Hochbruck2014, Kronbichler2016, Wang2019},
\begin{align}
\begin{aligned}
p_h^*&=\DGaverage{p_h}+\frac{\tau}{2} \DGjjump{\vec{u}_h},\\
\vec{u}_h^*&=\DGaverage{\vec{u}_h}+ \frac{\gamma}{2} \DGjjump{p_h}.\label{eq:fluxes}
\end{aligned}
\end{align}
The penalty parameters $\tau$ and $\gamma$ are derived using the Rankine--Hugoniot condition when solving for the Riemann solution~\cite{LeVeque2002, Hochbruck2014}, resulting in $\textstyle \tau=\rho c$ and~$\textstyle \gamma=\frac{1}{\rho c}$.
These penalty parameters are consistent in terms of a dimension analysis which demands~$\textstyle \tau\sim\rho c$ and~$\textstyle \gamma\sim\frac{1}{\rho c}$.
Element boundaries are either located inside the domain~$\pd\Omega^e_{\mathrm{inner}}$, at non-conforming boundaries~$\Gamma^e_{\mathrm{NCI}}$, or subject to BCs ($\pd\Omega_p^{\mathrm{D},e}$, $\pd\Omega_u^{\mathrm{D},e}$, or $\pd\Omega^{\mathrm{Y},e}$).
The explicit notation of the discretization at element boundaries reads
\begin{align}
\begin{aligned}
  \bigg(\frac{1}{\rho}\vec{w}_h\cdot\vec{n}, p_h^*\bigg)_{\pd\Omega^e}
  &=\left(\frac{1}{\rho}\vec{w}_h\cdot\vec{n}, p_h^*\right)_{\pd\Omega^e_{\mathrm{inner}}}
  +\left(\frac{1}{\rho}\vec{w}_h\cdot\vec{n}, g_p\right)_{\pd\Omega_p^{\mathrm{D},e}}
  +\left(\frac{1}{\rho}\vec{w}_h\cdot\vec{n}, p_h^-+\tau(\vec{u}_h^--\vec{g}_u)\right)_{\pd\Omega_u^{\mathrm{D},e}}
  \\
  &+\bigg(\frac{1}{\rho}\vec{w}_h\cdot\vec{n}, p_h^-+\tau\bigg (\vec{u}_h^-\cdot\vec{n}-\frac{Y}{\rho c}p_h^-\bigg)\bigg)_{\pd\Omega^{\mathrm{Y},e}}
+\left(\frac{1}{\rho}\vec{w}_h\cdot\vec{n}, p_{h,\mathrm{NCI}}^*(p_h^-,p_h^+,\vec{u}_h^-,\vec{u}_h^+))\right)_{\Gamma^e_{\mathrm{NCI}^-}},
\end{aligned} 
\end{align} 
and 
\begin{align}
\begin{aligned}
  \big(\rho c^2q_h\vec{n}, \vec{u}_h^*\big)_{\pd\Omega^e}
  &=\left(\rho c^2q_h\vec{n}, \vec{u}_h^*\right)_{\pd\Omega^e_{\mathrm{inner}}}
  +\left(\rho c^2q_h\vec{n}, \vec{u}_h^-+\gamma(p_h^--g_p)\right)_{\pd\Omega_p^{\mathrm{D},e}}
  +\left(\rho c^2q_h\vec{n}, \vec{g}_u\right)_{\pd\Omega_u^{\mathrm{D},e}}\\
  &+\left(cq_h, 2Yp_h^-\right)_{\pd\Omega^{\mathrm{Y},e}}
  +  \Big(\rho {c}^2q_h\vec{n},\vec{u}_{h,\mathrm{NCI}}^*(p_h^-,p_h^+,\vec{u}_h^-,\vec{u}_h^+)\Big)_{\Gamma^e_{\mathrm{NCI^-}}}.
\end{aligned} 
\end{align}
In this notation~$\Gamma^e_{\mathrm{NCI^-}}$ are the faces of the primary elements.
A more detailed description of how integrals are computed in the non-conforming case is provided in Section~\ref{sec:ptop} and \ref{sec:nitschemortar}.

To be able to consider different materials, we have to adapt the fluxes at the NCIs.
To this end, we use the LDG fluxes with special self-adapting upwinding parameters and penalty terms as an additional stabilization mechanism to increase the numerical diffusion~\cite{Hochbruck2014}  
\begin{align}
\begin{aligned}
p_{h,\mathrm{NCI}}^*&=p_h^--\frac{\tau^-}{\tau^-+\tau^+}\DGjump{p} +\frac{\tau^-\tau^+}{\tau^-+\tau^+} \DGjjump{\vec{u}_h},
\\
\vec{u}_{h,\mathrm{NCI}}^*&=\vec{u}_h^--\frac{\gamma^-}{\gamma^-+\gamma^+}\DGjump{\vec{u}_h}+ \frac{\gamma^-\gamma^+}{\gamma^-+\gamma^+}\DGjjump{p_h}.
\end{aligned}
\end{align}
We can see, that the fluxes simplify to the Lax--Friedrichs fluxes of~\eqref{eq:fluxes} for homogenous materials, i.e. $c^-\rho^-$=$c^+\rho^+$ and therefore, $\gamma^-=\gamma^+$ and $\tau^-=\tau^+$.

\subsection{Point-to-Point Interpolation\label{sec:ptop}}
\begin{figure}[t]
  \centering
  \includegraphics[width=.4\textwidth]{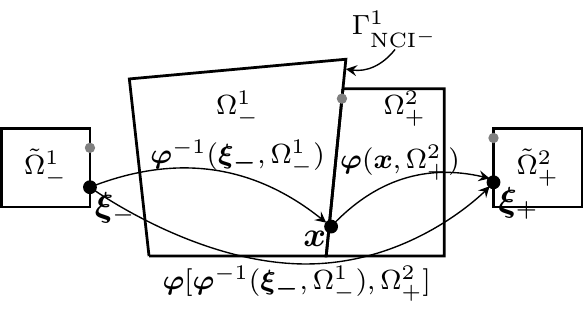}
  \caption{Point-to-point interpolation: Shown is the mapping of a exemplary quadrature point~$\vec{\xi}_-$ (associated to the primary element~$\Omega^1_{-}$) to the the non-conformingly connected secondary element~$\Omega^2_{+}$. Adapted from~\cite{Heinz2022}.\label{fig:mapping}}
\end{figure}

Non-conformity can be easily handled by the evaluation of all quantities in consistent quadrature points~\cite{Hermann2010, Laughton2021};
i.e., we have to evaluate fluxes at the same point in the physical space.
The primary element dictates the used quadrature points; see Figure~\ref{fig:mapping}. 
For conforming DG this leads to the same quadrature points in the reference space $\vec{\xi}_-=\vec{\xi}_+$.
However, if non-conformities in the mesh are present, quadrature points that correspond to the same coordinate in the physical space differ, and we have to find quadrature points on the secondary element as
\begin{align}
  \vec{\xi}^+=\vec{\varphi}(\vec{x},\Omega^e_+).
\end{align} 
Therefore, we can explicitly state that an arbitrary physical flux~$\mathcal{F}_h^*$ is a function of arbitrary quantities~$b$, evaluated at the same physical coordinates~$\vec{x}=\varphi^{-1}(\vec{\xi}^-,\Omega_-^e)$ (provided by the integration rule of the corresponding primary element face)
\begin{align}
  \begin{aligned}
    \mathcal{F}_h^*(b^-,b^+)=\mathcal{F}_h^*\left(b^-(\vec{\xi}^-),b^+(\vec{\xi}^+)\right).
\end{aligned}
\end{align}
The integration over a non-conforming face of a primary element subsequently reads
\begin{align}
\begin{aligned}
  \int_{\pd\Omega} \mathcal{F}^*(b^-,b^+)\ \dd\Gamma_{\mathrm{NCI^-}}
  \approx\sum_{q=1}^{n_q} w_q \mathcal{F}^*_h(b^-(\vec{\xi}^-),b^+(\vec{\xi}^+))|J_q|.
\end{aligned}
\end{align}
$w_q$ are the weights of the Gauss quadrature, and the Jacobi determinants~$|J_q|$ in quadrature points~$q$ correspond to the primary element face.
\subsection{Nitsche-Type Mortaring\label{sec:nitschemortar}}
\begin{figure}[t]
  \centering
  \includegraphics[width=.4\textwidth]{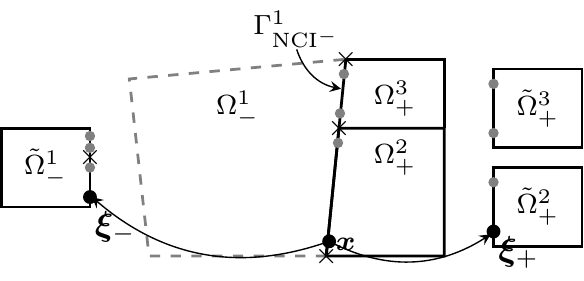}
  \caption{Nitsche-type mortaring: In contrast to point-to-point interpolation, quadrature points are not dictated by the primary element.
  Instead, the mortars between elements are computed, and quadrature points, weights, and Jacobians correspond to the mortar patches.
  This way, there is no aliasing, and values form a smooth representation in each quadrature.
  Intersections are computed between the face of the primary element and the secondary volume elements to ensure the method works for overlaps without modification.\label{fig:mapping_mortaring}}
\end{figure}
In contrast to point-to-point interpolation, Nitsche-type mortaring computes the integration over primary element faces at NCIs as the sum of collected quadrature rules on the mortars, see Figure~\ref{fig:mapping_mortaring}.
Thus the integral computes as
\begin{align}
\begin{aligned}
  \int_{\pd\Omega} \mathcal{F}^*(b^-,b^+)\ \dd\Gamma_{\mathrm{NCI^-}}
  \approx\sum_{m=1}^{n_m} \sum_{q=1}^{n_q} w_q \mathcal{F}^*_h(b^-(\vec{\xi}^-),b^+(\vec{\xi}^+))|J_q^m|.\label{eq:mortarint}
\end{aligned}
\end{align}
In~\eqref{eq:mortarint},~$n_m$ is the number of found intersections.
The Jacobi determinants~$|J_q^m|$ in quadrature point~$q$ is now determined on mortar~$m$.
This way, the integral on the NCI is computed exactly (if enough quadrature points~$n_q$ are chosen) without aliasing.
Constructing the mortars is more challenging to implement and reduces performance~\cite{Laughton2021}.
In our formulation, intersections are computed between the face of a primary element and the secondary volume elements, cf. Figure~\ref{fig:mapping_mortaring};
in 2D simulations, intersections are computed between a quadrilateral and a line.
This way, we can seamlessly handle overlaps between elements (see Section~\ref{sec:overlap}).

\section{Remarks on Implementation\label{sec:imple}}
Our implementations will be freely available as a part of the open source software project \texttt{EXADG}~\cite{Arndt2020} and the software library \texttt{deal.II}~\cite{Arndt2022}.

In the case of point-to-point interpolation, we are collecting all quadrature points on the NCIs mapped to the physical space.
After that, we perform a global search based on distributed bounding boxes for secondary elements that hold these integration points and store the corresponding quadrature points in the reference space.
In each time step, we evaluate pressure and velocities on the secondary elements and use the results to compute the fluxes over the NCIs.
If a quadrature point is found on multiple elements, we use the average value in the computation of the fluxes.
This approach works on curved elements without further ado.

In the case of Nitsche-type mortaring, we first create the mortars between the primary and secondary elements.
We are computing the $d-1$ dimensional intersections between $d$ and $d-1$ geometric entities.
The $d-1$ element is a face of a primary element.
This way, quadrature rules are defined on the primary element faces, independent if elements overlap or not.
Thus, our implementation is the same in case of element overlaps or standard NCIs.
In our implementations, we extract the vertices of the non-conforming faces of the primary elements and all vertices of possibly connected secondary elements.
We then use \texttt{CGAL}~\cite{CGAL2022} to compute the inter-dimensional intersections and eventually create mapped quadrature rules on each mortar patch.
The rest of the implementation is nearly the same as for the point-to-point interpolation:
Additionally to the quadrature points, we store the Jacobi determinants at the quadrature points.
Normal vectors are not stored; instead, we use the negative normal vectors of the primary element during flux evaluation.
We evaluate pressure and velocity on the primary and secondary elements in each time step and use the stored Jacobi determinants to compute and test the fluxes over the NCIs.
This approach limits us to non-curved elements at NCIs.
Since \texttt{CGAL} is working with triangular/tetrahedral elements, mortar patches are always triangular in the 3D case (even if the patch could be rectangular).
Note that the number of created mortars and thus the number of quadrature points highly depends on the element configuration.

\section{Numerical Results\label{sec:numeric}}

\begin{figure}[t] 
  \centering  
  \includegraphics[width=0.35\textwidth]{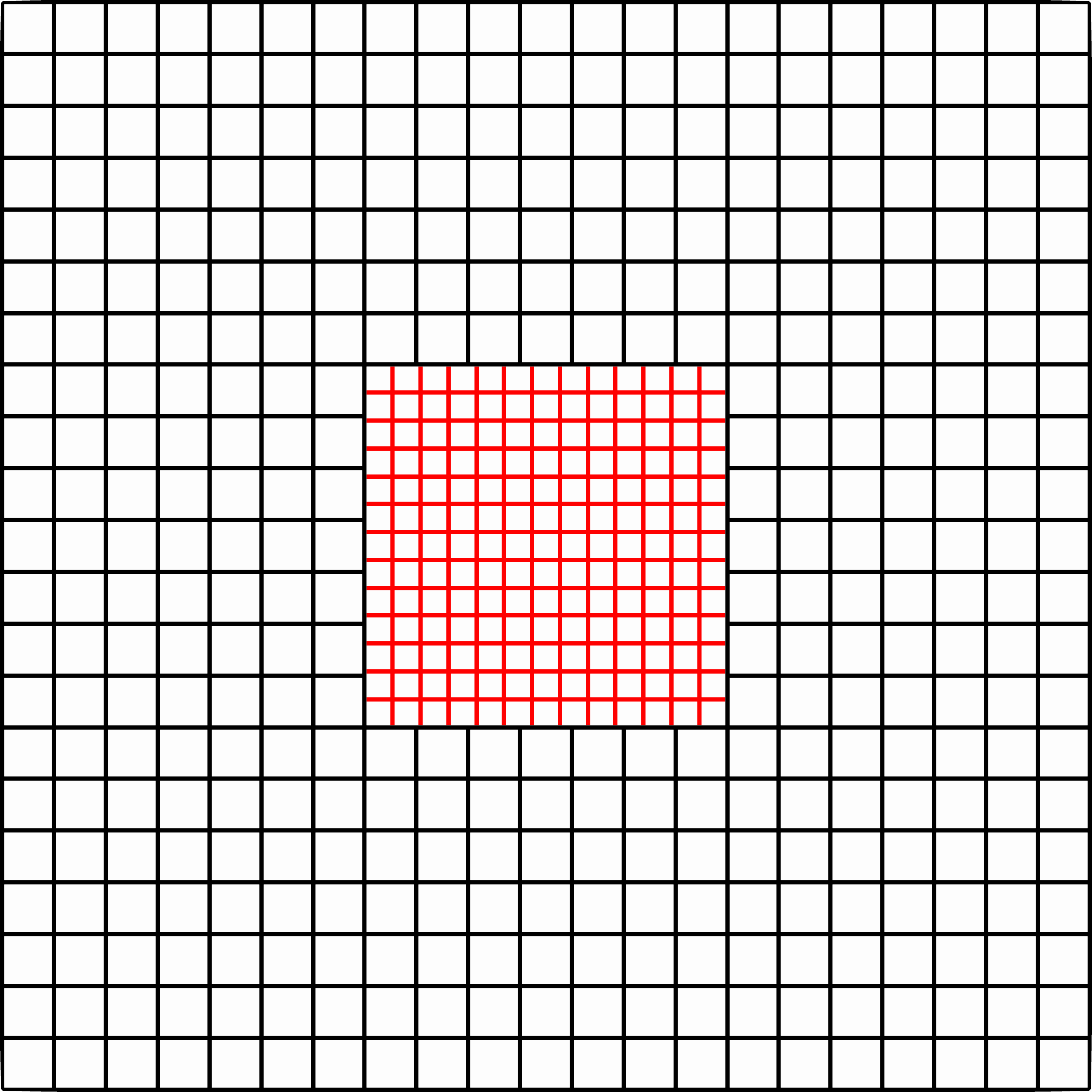}
  \caption{Computational mesh consisting of two mesh regions which are connected via a non-conforming interface.\label{fig:mesh}}
\end{figure}  
This section provides numerical results using point-to-point interpolation and Nitsche-type mortaring.
First, we show instabilities related to non-smooth representations of values at NCIs for the point-to-point interpolation method.
We show that these instabilities do not occur if we use Nitsche-type mortaring instead.
For Nitsche-type mortaring, we provide in-depth convergence results for different mesh regions, quantify the error introduced at the NCI, and provide results for a test case with heterogeneous material.
Additionally, we show that the method also works if elements are overlapping.
From now on, all spatial values are given in$\unit{m}$.
\subsection{Vibrating Membrane} 
To be able to compute errors exactly, we use the test case of a vibrating membrane which has also been used, amongst others, in~\cite{Nguyen2011, Schoeder2019}.
For this test case, the analytical solution at each time $t$ is known, and for a two-dimensional domain, it reads for the pressure
\begin{align}
  p =\cos(M \sqrt{2} \pi t) 
  \sin(M \pi x) \sin(M \pi y),\label{eq:pana}
\end{align}
and for the acoustic particle velocity
\begin{align} 
\vec{u}=\frac{-\sin(M  \sqrt{2} \pi t)}{\sqrt{2}}
  \begin{pmatrix}
    \cos(M \pi x) \sin(M \pi y) \\
    \sin(M \pi x) \cos(M \pi y)
\end{pmatrix},
\end{align}
assuming no acoustic loads~$F=0\unit{kg\cdot m^{-3}\cdot s^{-1}}$, as well as~$\rho=1\unit{kg\cdot m^{-3}}$ and~$c=1\unit{m\cdot s^{-1}}$.
Our simulations' computational domain $\Omega$ consists of two mesh regions connected via NCIs, as shown in Figure~\ref{fig:mesh}.
The outer region $\Omega_{\mathrm{o}}$ is a rectangular domain with a rectangular hole in which the inner region $\Omega_{\mathrm{i}}$ is embedded.
Thus, $\Omega = \Omega_{\mathrm{o}}\cup\Omega_{\mathrm{i}}$ and within the following $\Omega_{\mathrm{o}}=(0,0)\times(0.1,0.1)\setminus \Omega_{\mathrm{i}}$ and $\Omega_{\mathrm{i}}=\left(\nicefrac{1}{30},\nicefrac{1}{30}\right)\times\left(\nicefrac{2}{30},\nicefrac{2}{30}\right)$. 
We use $M=120$ modes, which leads to $p=0\unit{Pa}$ at the boundaries of the computational domain $\Omega$ and we apply homogenous pressure DBCs $g_p=0\unit{Pa}$.
All computations use the low storage Runge--Kutta version RKC84~\cite{Toulorge2012}.
The CFL condition
\begin{align}
  \Delta t=\frac{\mathrm{Cr}}{k^{1.5}}\frac{h_{\min}}{c_{\max}},\label{eq:clf}
\end{align}
gives the required time step size $\Delta t$ needed for a stable temporal discretization.
We use the minimum edge length~$h_{\min}$ of all existing elements, the polynomial degree $k$, and the largest value of the speed of sound $c_{\max}$ for its computation.
To account for the different spacing between Legendre--Gauss--Lobatto (used as quadrature points) we are using the superscript~$\textstyle 1.5$ as proposed in~\cite{Fehn2018}.
For a detailed discussion on the CFL condition for explicit Runge--Kutta methods, we refer to~\cite{Toulorge2011}.
From this point forward, all computations use time step sizes computed by the CFL condition with a Courant number $\mathrm{Cr}=0.2$. 

\subsubsection{Instabilities\label{sec:instab}}
\begin{figure*}[t]
  \centering
  \subfloat[][Point-to-point interpolation.]{\includegraphics[width=.5\textwidth]{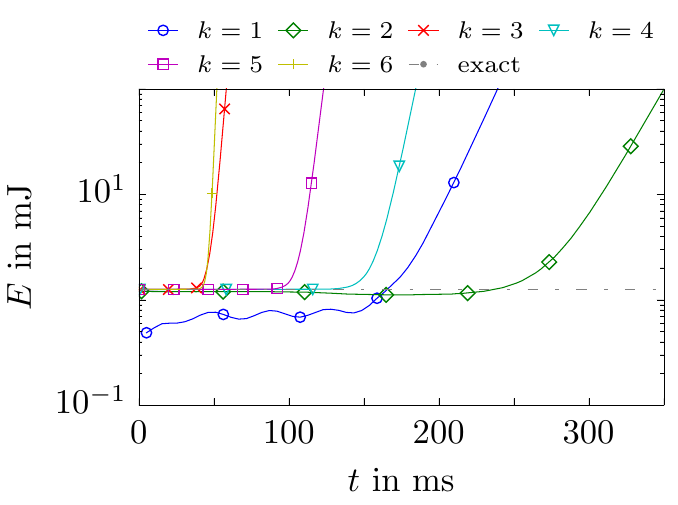}\label{fig:instab}}\hfill
  \subfloat[][Nitsche-type mortaring.]{\includegraphics[width=.5\textwidth]{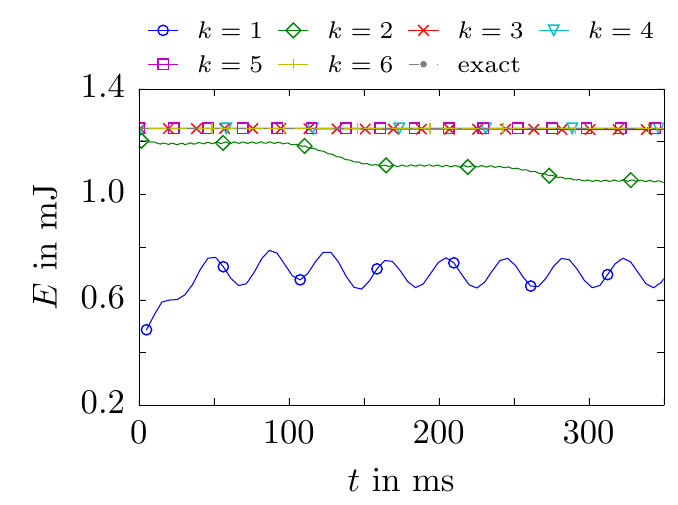}\label{fig:overint}}
  \caption{Sound energy over time computed for orders $k=1$ to $k=6$. The given setup for the vibrating membrane test case is perfectly energy conservative and thus energy has to be constant over time. For point-to-point interpolation, instabilities form after some time while the simulation stays stable for exact integration with Nitsche-type mortaring.\label{fig:label}} 
\end{figure*}
The test case perfectly conserves the sound energy  
\begin{align}
E=\int_{\Omega}\left(\frac{p^2}{2\rho c^2} + \frac{\rho (\vec{u}\cdot\vec{u})}{2}\right) \dd \Omega. 
\end{align}
Since the analytic solution exists, we can compute the exact sound energy contained in the system as $E_{\mathrm{exact}}=1.25\unit{mJ}$.
The mesh (cf. Figure~\ref{fig:mesh}) has element edge lengths of $h_{\Omega_\mathrm{i}}=\nicefrac{1}{30\cdot 13}$ on the NCI for the inner domain and $h_{\Omega_\mathrm{i}}=\nicefrac{1}{30\cdot 7}$ on the NCI for the outer domain.
Figure~\ref{fig:instab} shows the sound energy in the system over time for point-to-point interpolation.
After a certain time, instabilities manifest as an non-physical rapid increase of sound energy.

Obviously, the approach suffers from aliasing; 
the integration of the primary elements only includes information from each connected element if quadrature points are found in every element.
One can regard this as a Dirichlet-Dirichlet approach, where the values are picked from the secondary domain instead of, e.g., an analytical function.
This reasoning does not explain the observed instabilities.

The difference in the applied Dirichlet boundary values is that in the case of an analytical function, the boundary values form a smooth representation of the solution.
In the case of using values from the secondary domain, boundary values are not necessarily smooth.
If quadrature points are located in different elements, we might face jumps in the solution representation.
While these jumps are assumed to be less distinctive in the case of continuous FE methods, the nature of the DG method intensifies this issue.
Nevertheless, the same also happens in the case of continuous FE methods, in the case where whole secondary elements are not sampled by any quadrature point of the primary element.
The jumps between Dirichlet values introduce spurious oscillations that eventually lead to unstable simulations.
To quantify that this is indeed the source of instabilities, we tested to interpolate solution values between domains into the DoFs.
This way, there are no jumps between quadrature points since the shape function of the primal element enforces a smooth representation of the values.
Nevertheless, we observed instabilities in the case of high polynomial degrees.
These instabilities are related to Gibbs' phenomena.
Significant differences between DoFs and the high-order shape functions lead to ringing modes, typically observed in shock capturing.
Applying techniques to interpolated values that are usually used in shock capturing, such as modal filtering~\cite{Hesthaven2008}, lead to stable schemes.
However, this also leads to a drop in the obtained spatial convergence rates; 
therefore, we are not further discussing those approaches.

According to the previous discussion, the phenomenon is not expected in the case of Nitsche-type mortaring.
Since the integrals are evaluated on intersections between elements, the representation of values is forced to be continuous in each quadrature rule.
The results are depicted in Figure~\ref{fig:overint}; we do not observe instabilities.
From this point forward, we will only consider the Nitsche-type mortaring approach.

Note that the main focus of the test case is to show instabilities for any polynomial degree.
However, the mesh resolution for low polynomial degrees $k$ is chosen too poor for a good approximation of the primal variables.
According to~\cite{Ainsworth2004}, the element size $h$ has to be chosen such that
\begin{align}
k+\frac{1}{2}>\frac{\omega h}{2}+C(\omega h)^{1/3}.
\label{eq:resolution}
\end{align}
In~\eqref{eq:resolution}, $C$ is a constant that can be chosen unity in practice~\cite{Ainsworth2004} and $\omega$ is the wave-number.
For the vibrating membrane test case $\omega=2\pi M$.
Using the maximum element size in used triangulation, we obtain $\frac{\omega h}{2}+(\omega h)^{1/3}\approx 3.44$ for current investigations.
Thus, we fulfill~\eqref{eq:resolution} with polynomial degrees $k>2$, for $k=2$ we are slightly off, and for $k=1$ we have a substantial deviation.
Consequently, we can see constant sound energy over time for orders greater than $k=2$.
For $k=2$, we see a slight drop in energy due to numerical dispersion originating from the slightly too coarse resolution.
For $k=1$, the resolution is so poor that the sound energy can not be resolved from the beginning, and we observe non-physical oscillations.
These results are in accordance to~\cite{Ainsworth2004}, in which it is reported that results might even become qualitatively incorrect for insufficient resolutions.
Within the next sections the mesh resolutions are chosen such that~\eqref{eq:resolution} is fulfilled for $k=1$ which successfully removes any non-physical oscillations.

\subsubsection{Convergence Results}
\begin{figure*}[p] 
  \centering
  \subfloat[][Global domain~$\mathcal{R}$~is~$\Omega$.]{\includegraphics[width=\textwidth]{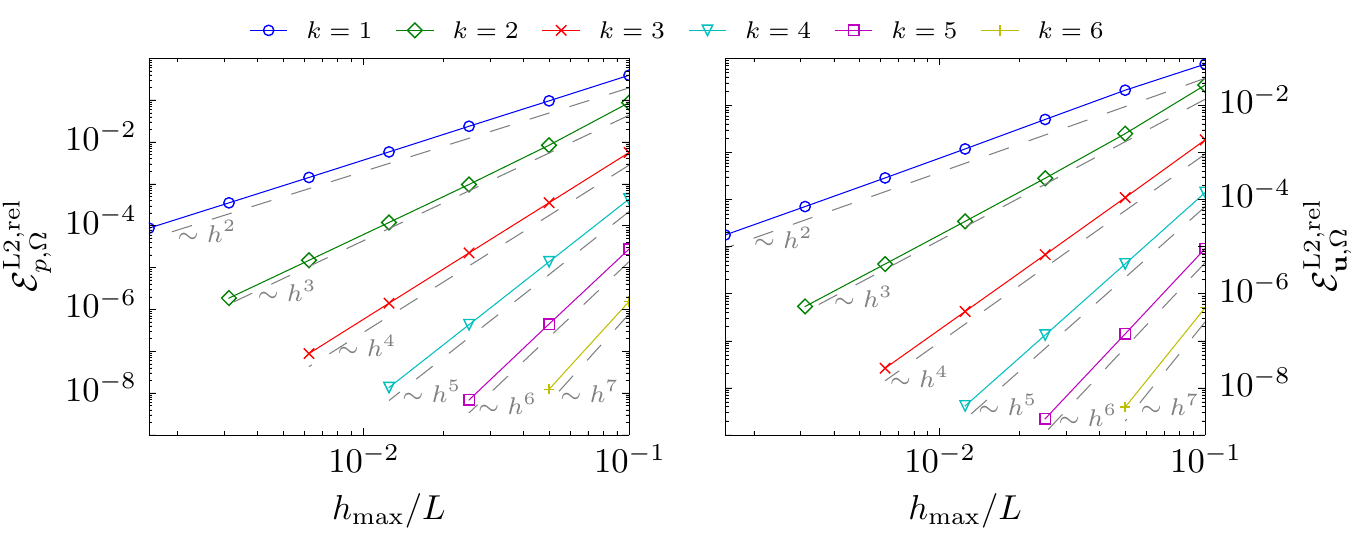}\label{fig:convergence_global}}
  
  \subfloat[][Outer domain~$\mathcal{R}$~is~$\Omega_{\mathrm{o}}$.]{\includegraphics[width=\textwidth]{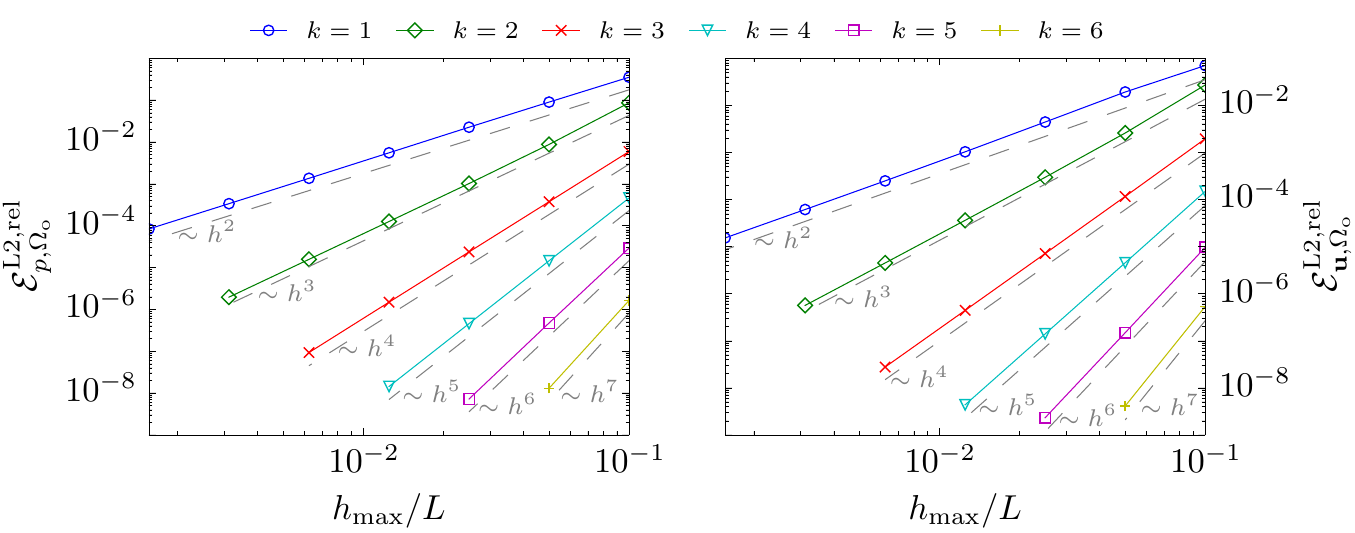}\label{fig:convergence_outer}}
  
  \subfloat[][Inner domain~$\mathcal{R}$~is~$\Omega_{\mathrm{i}}$.]{\includegraphics[width=\textwidth]{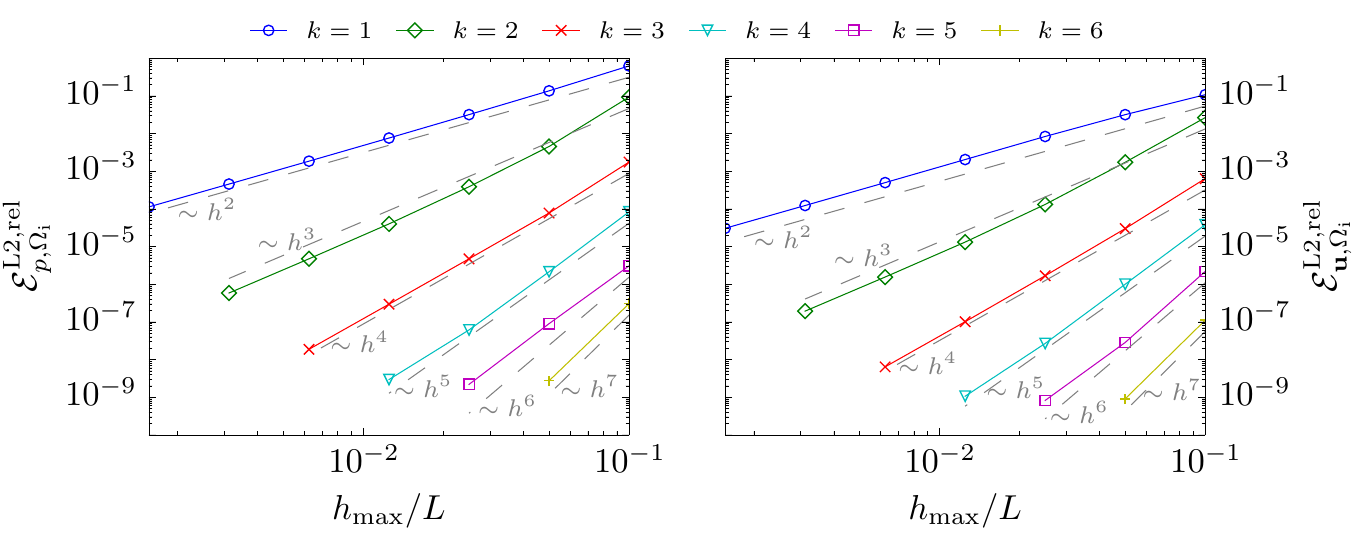}\label{fig:convergence_inner}}
  \caption{Spatial convergence study for the vibrating membrane testcase with $M=30$, defined on the rectangular domain of Figure~\ref{fig:mesh} (with $h_{\mathrm{i,initial}}=\nicefrac{1}{30\cdot 3}$, and $h_{\mathrm{o,initial}} =\nicefrac{1}{30\cdot 2}$) using Nitsche-type mortaring: Shown are the relative~$\mathrm{L2}$ errors for pressure~$\mathrm{\varepsilon}^{\mathrm{L2,rel}}_{p,\mathcal{R}}$ and velocity~$\mathrm{\varepsilon}^{\mathrm{L2,rel}}_{\vec{u},\mathcal{R}}$ on different domains~$\mathcal{R}$. The domain $\mathcal{R}$ might be the global domain $\Omega$, the inner domain $\Omega_{\mathrm{i}}$, or the outer domain $\Omega_{\mathrm{o}}$.\label{fig:convergence}}
\end{figure*} 

We use the setup of the previous section with $M=30$ modes.
However, we alter the mesh sizes compared to Figure~\ref{fig:mesh}.
The elements in the outer domain have initial edge lengths of $h_{\mathrm{o,initial}} =\nicefrac{1}{30\cdot 2}$;
in the inner domain initial element edge lengths are $h_{\mathrm{i,initial}}=\nicefrac{1}{30\cdot 3}$.
We compute the relative $L2$ error for the pressure $\mathrm{\varepsilon}^{\mathrm{L2,rel}}_{p,\mathcal{R}}$ on region $\mathcal{R}$ after $1\unit{s}$ for different mesh refinements
\begin{align}
\mathrm{\varepsilon}^{\mathrm{L2,rel}}_{p,\mathcal{R}}=\frac{\sqrt{\int_{\mathcal{R}}  (p_h-p_{\mathrm{ana}})^2\ \dd \Omega}}{\sqrt{\int_{\mathcal{R}}  p_{\mathrm{ana}}^2\ \dd \Omega}},
\end{align}
with the analytical solution of the pressure $p_{\mathrm{ana}}$, see Eq.~\eqref{eq:pana}.
The velocity error $\mathrm{\varepsilon}^{\mathrm{L2,rel}}_{\vec{u},\mathcal{R}}$ is computed accordingly.
Regions are either the global region $\Omega$, the inner region $\Omega_{\mathrm{i}}$, or the outer region  $\Omega_{\mathrm{o}}$.
The mesh refinement is realized by replacing each quadrilateral cell by four children cells, and the corresponding edge lengths $h$ at refinement level $r$ compute as
\begin{align}
h=\frac{h_{\mathrm{initial}}}{2^{r}}.
\end{align}

We observe optimal convergence rates of order $k+1$ in space on the global domain $\Omega$, see Figure~\ref{fig:convergence_global}.
The outer domain has a coarser spatial discretization and dominates the errors on the global domain.
Therefore, it is not surprising that the errors on the outer domain (Figure~\ref{fig:convergence_outer}) behave similar to the ones on the global domain (Figure~\ref{fig:convergence_global}).
The inner domain has a finer spatial discretization;
thus the errors obtained in the inner domain might be shadowed by the errors obtained in the outer domain.
However, errors propagate from the outer domain to the inner domain.
We also observe optimal convergence rates computing the errors on the inner domain (see Figure~\ref{fig:convergence_inner}).
Errors obtained in the inner and outer domain are similar.
Thus errors from the outer domain entirely propagated to the inner domain after $1\unit{s}$.
Therefore, in practical applications, one should aim to choose mesh sizes that yield approximately the same errors in each domain.
In conclusion, we obtain optimal convergence rates in all regions of the non-conforming mesh and can confidently apply the proposed method, keeping in mind that a jump in element sizes has to be justified, e.g., due to different materials.

\subsubsection{Embedding of Circular Domain\label{sec:overlap}}
\begin{figure*}[!t]
  \centering
  \subfloat[][Curved interface]{\includegraphics[width=.32\textwidth]{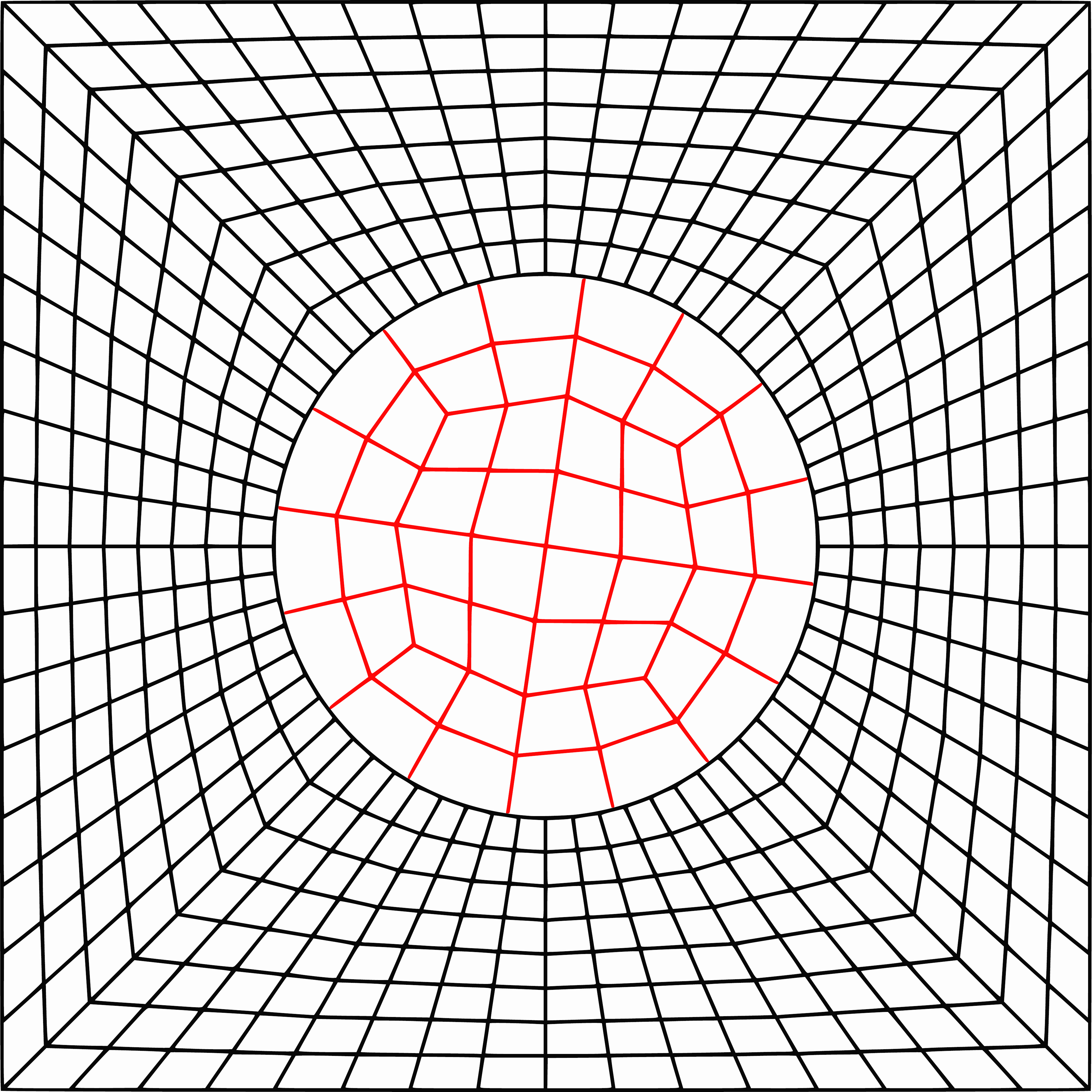}\label{fig:curved}}\hfill
  \subfloat[][Small overlap]{\includegraphics[width=.32\textwidth]{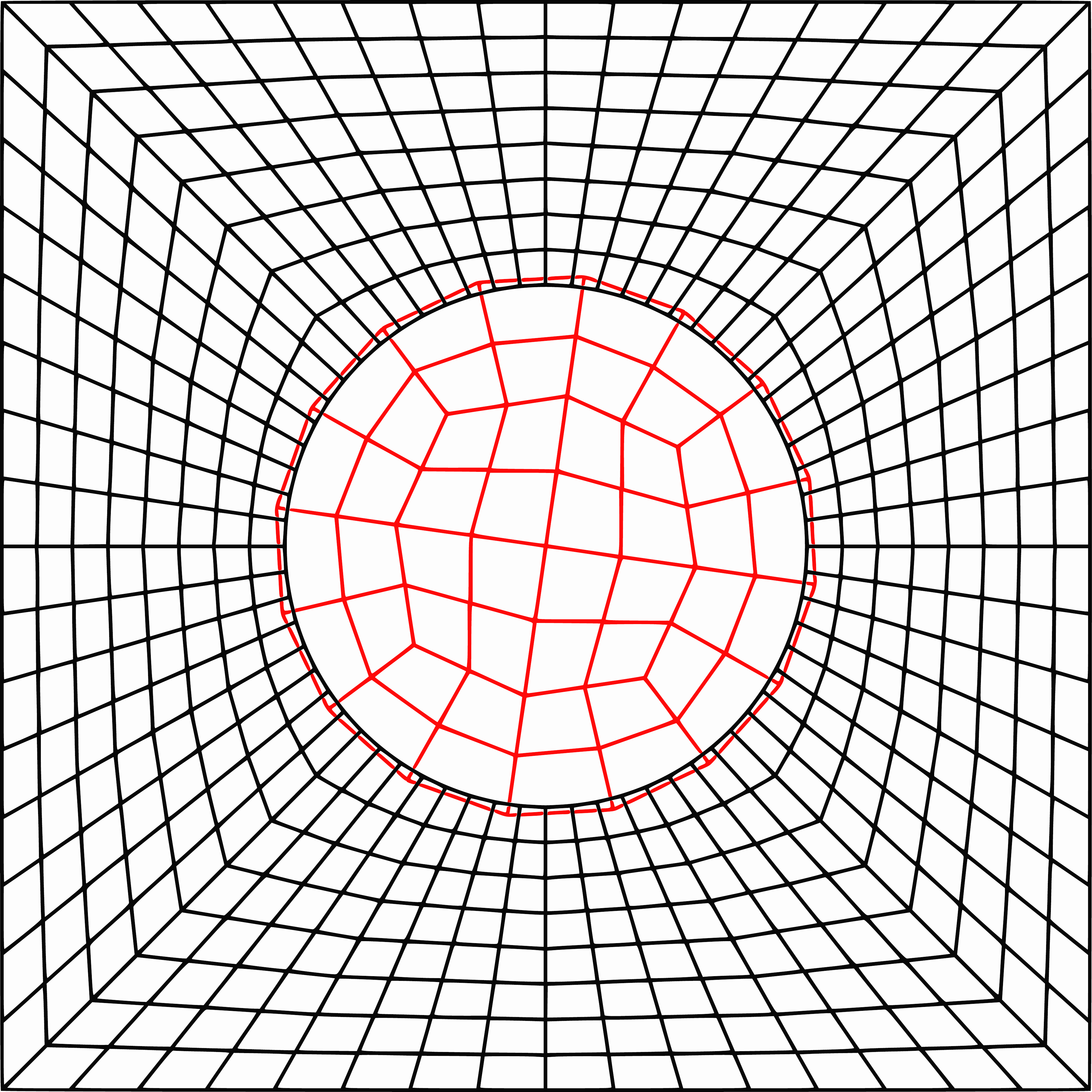}\label{fig:overlap}}\hfill
  \subfloat[][Overset mesh]{\includegraphics[width=.32\textwidth]{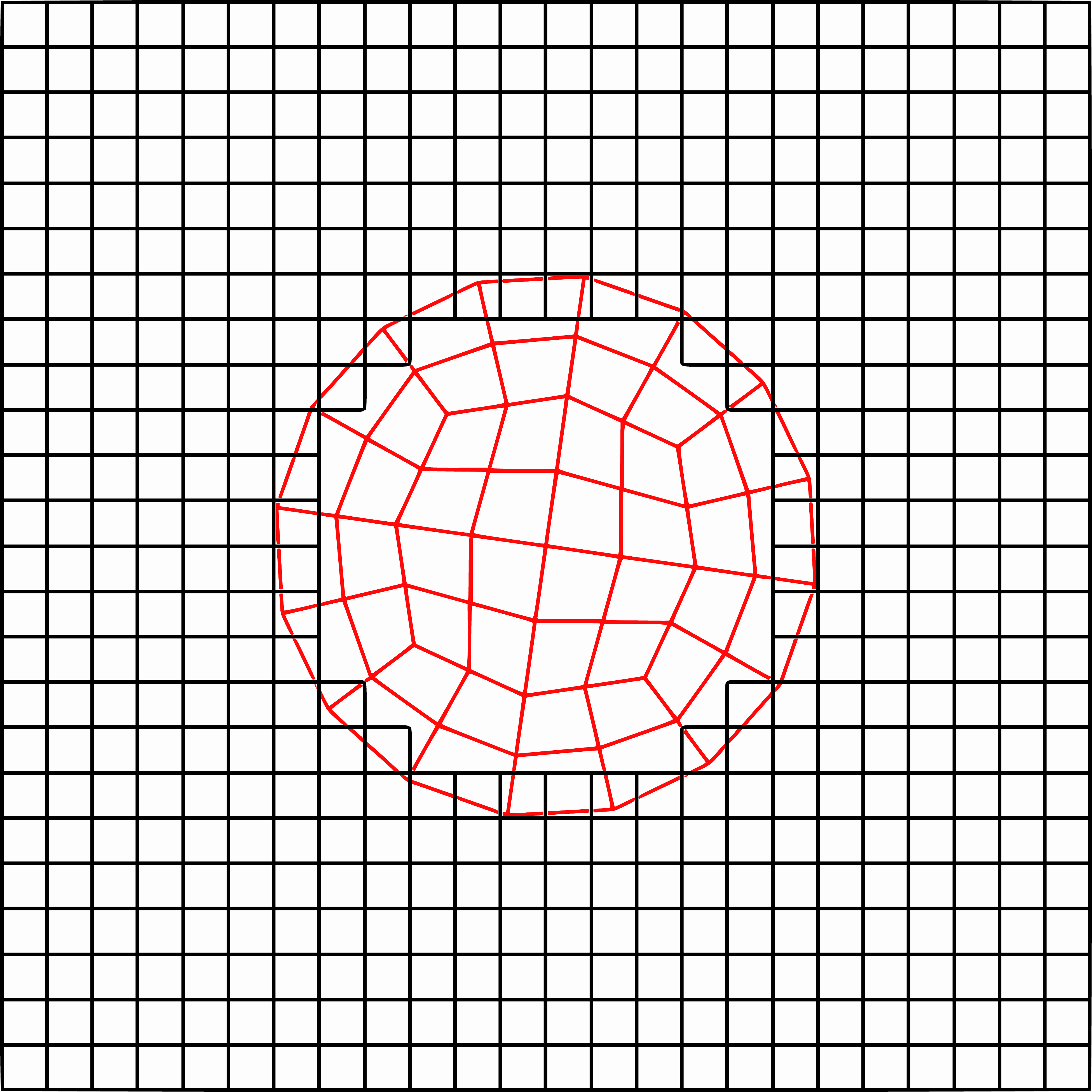}\label{fig:overset}}\hfill
  \caption{Three different meshes that have a circular mesh embedded in a rectangular mesh.\label{fig:mesh_overlap}} 
\end{figure*}

\begin{table*}[b]
    \centering
    \caption[]{Relative~$L^2$ errors~$\varepsilon^{\mathrm{L2,rel}}_{\Omega}=\varepsilon^{\mathrm{L2,rel}}_{p,\Omega}+\varepsilon^{\mathrm{L2,rel}}_{\vec{u},\Omega}$ for different polynomial degrees $k$ computed on the meshes depicted in Figure~\ref{fig:mesh_overlap}.\label{tab:errorsoverlap}}
    \begin{tabular*}{500pt}{@{\extracolsep\fill}lcc@{\extracolsep\fill}cccc@{\extracolsep\fill}}
        \toprule
       & \multicolumn{2}{@{}c@{}}{\textbf{Curved interface}} & \multicolumn{2}{@{}c@{}}{\textbf{Overlap}} & \multicolumn{2}{@{}c@{}}{\textbf{Overset}} \\
        \cmidrule{2-3} \cmidrule{4-5} \cmidrule{6-7}
        $k$&\textbf{DoFs} & $\varepsilon^{\mathrm{L2,rel}}_{\Omega}$ &\textbf{DoFs} & $\varepsilon^{\mathrm{L2,rel}}_{\Omega}$&\textbf{DoFs} & $\varepsilon^{\mathrm{L2,rel}}_{\Omega}$ \\ 
        \midrule 
        1 &6,720 & $2.915\cdot 10^{-2}$ & 6,720  & $2.779\cdot 10^{-2}$ & 6,432  & $2.053\cdot 10^{-2}$  \\
        2 &15,120& $8.303\cdot 10^{-4}$ & 15,120 & $5.444\cdot 10^{-4}$ & 14,472 & $4.641\cdot 10^{-4}$  \\
        3 &26,880& $4.304\cdot 10^{-5}$ & 26,880 & $8.536\cdot 10^{-6}$ & 25,728 & $7.428\cdot 10^{-6}$  \\
        4 &42,000& $1.666\cdot 10^{-6}$ & 42,000 & $2.220\cdot 10^{-7}$ & 40,200 & $2.200\cdot 10^{-7}$  \\
        5 &60,480& $6.974\cdot 10^{-8}$ & 60,480 & $2.328\cdot 10^{-9}$ & 57,888 & $2.367\cdot 10^{-9}$  \\
        6 &82,320& $2.650\cdot 10^{-9}$ & 82,320 & $1.286\cdot 10^{-9}$ & 78,792 & $1.096\cdot 10^{-9}$  \\
        \bottomrule 
    \end{tabular*} 
  \end{table*}
   
Being able to handle overlaps has two useful properties.
Mesh generation gets more straightforward, and rotating interfaces can be handled without the need for curved elements.
We provide results for three different grids, depicted in Figure~\ref{fig:mesh_overlap}, that are prototypical in the context of rotating interfaces.
The rectangular domain spans~$\Omega_{\mathrm{o}}=[0,0.1]^2$ and the circular domain $\Omega_{\mathrm{i}}$ has a radius of $0.5$.

In this particular case, it is easily possible to manually compute quadrature rules on the curved intersections since the NCI is a circle (cf. Figure~\ref{fig:curved}).
Note that this does not work for arbitrary shapes in our implementations since we rely on \texttt{CGAL} to compute the intersections.
Nevertheless, this approach becomes relevant for large-scale computations with sliding interfaces since the computational cost to create mortars is heavily reduced, cf.~\cite{Duerrwaechter2021}.
For the version with a slight overlap (cf. Figure~\ref{fig:overlap}), the radius of the hole is slightly smaller than the radius of the circular domain, i.e., $0.5- 2 \cdot 10^{-3}$.
Using overset meshes (cf. Figure~\ref{fig:overset}) is particularly helpful in generating structured meshes in regions connected to complex geometries.

All meshes have similar numbers of DoFs.
Note that the methodology works for arbitrary overlaps.
However, the same physical fields are computed in the overlap; 
thus, redundant work is done if the overlap exceeds one element.
Table~\ref{tab:errorsoverlap} shows the errors obtained after~$1\unit{s}$ for the vibrating membrane test case with~$M=5$ modes.
In this case, we apply inhomogeneous pressure DBCs with~$g_p$ obtained from the analytical solution.
We can see that the errors are in the same order of magnitude for the overlapping and overset mesh.
Even though we used fewer DoFs in the overset mesh, we can see slightly better errors, with an outlier at polynomial degree $k=5$.
This relates to the element distortions in the overlapping case.
The curved interface setup produces more significant errors than the overlapping setup, the most distinct deviations are for polynomial degree $k=3$ and $k=4$.
This is not expected and needs further investigation before application to sliding rotating interfaces.
One possible explanation is that round-off errors are introduced while computing the curved intersections.

Overall, we conclude that our methodology works as expected if elements overlap.

\FloatBarrier
\subsection{Application}
\begin{figure}[!h]
  \centering
  \includegraphics[width=.35\textwidth]{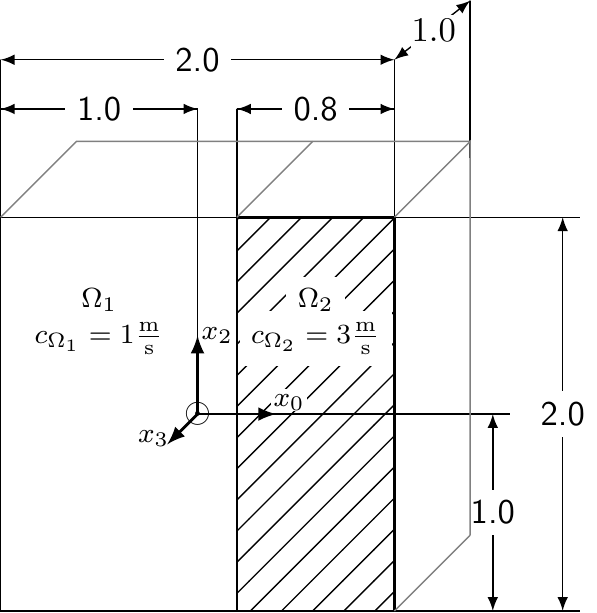}
  \caption{Application: Domain with heterogeneous fluids clipped in $x_1$-$x_2$ plane. The speed of sound in both fluids differs, while the density is $\rho_{\Omega_1}=\rho_{\Omega_2}=1\unit{kg\cdot m^{-3}}=\mathrm{const}$. A pressure pulse~$\odot$ is located in the center of the domain as an initial condition.\label{fig:appl_setup} }
\end{figure}
As pointed out, NCIs are especially desirable if different spatial resolutions are required.
Imagine two fluids with different speeds of sound $c$.
We need different element sizes to resolve the acoustic pressure up to a specific frequency.
We use the test case with heterogeneous acoustic material, also simulated by~\cite{Bangerth2010,Koecher2014,Perugia2020}.
We adapt the computational domain to show that our implementations work in the 3D case.
A wave travels over the interface between two materials.
At the interface, the wave is partially transmitted and partially reflected, and an additional wavefront emerges due to the Huygens--Fresnel principle.
A sectional view of the setup for this test case is depicted in Figure~\ref{fig:appl_setup}.
The domain $\Omega=\Omega_1\cup\Omega_2$ spans from $\Omega=(-1,-1,-1)\times(1,1,1)$.
In the left part of the domain the speed of sound is $c_{\Omega_1}=1\unit{m\cdot s^{-1}}$ while it is $c_{\Omega_2}=3\unit{m\cdot s^{-1}}$ in the right part.
The density of both fluids is $\rho_{\Omega_1}=\rho_{\Omega_2}=1\unit{kg\cdot m^{-3}}$.
As an initial condition, a pressure pulse is chosen
\begin{align}
p(t=0)&=\exp{-10^{4}\ \vec{x}\cdot\vec{x}},\\
\vec{u}(t=0)&=\vec{0}.
\end{align}
The test case is subject to homogenous pressure BCs.

\begin{figure*}[b] 
  \centering
  \subfloat[][$x_1$-$x_2$ plane, $t=0.1\unit{s}$.]{\includegraphics[width=0.48\textwidth]{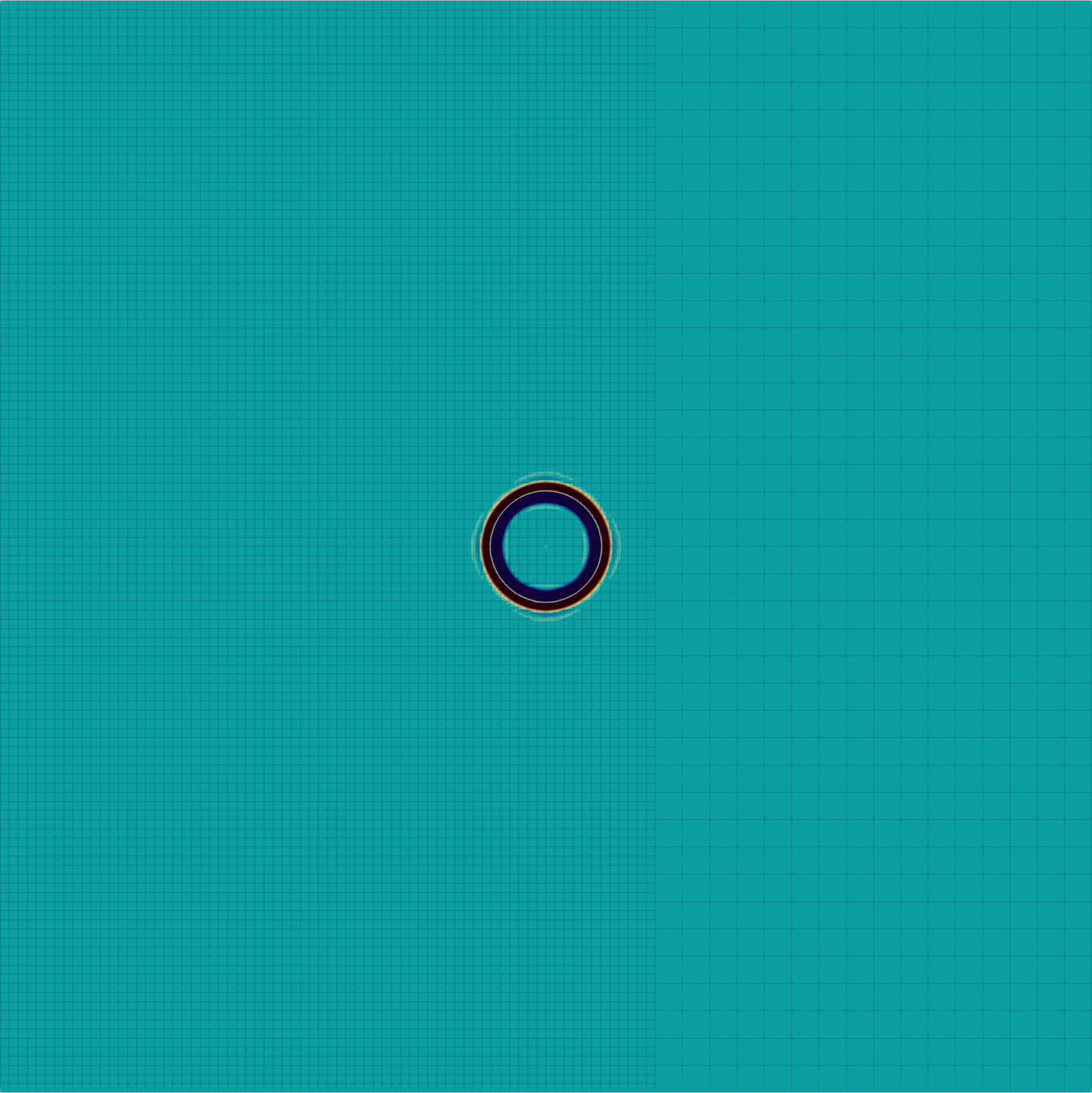}\label{fig:appl1}}\hfill
  \subfloat[][$x_1$-$x_2$ plane, $t=0.2\unit{s}$.]{\includegraphics[width=0.48\textwidth]{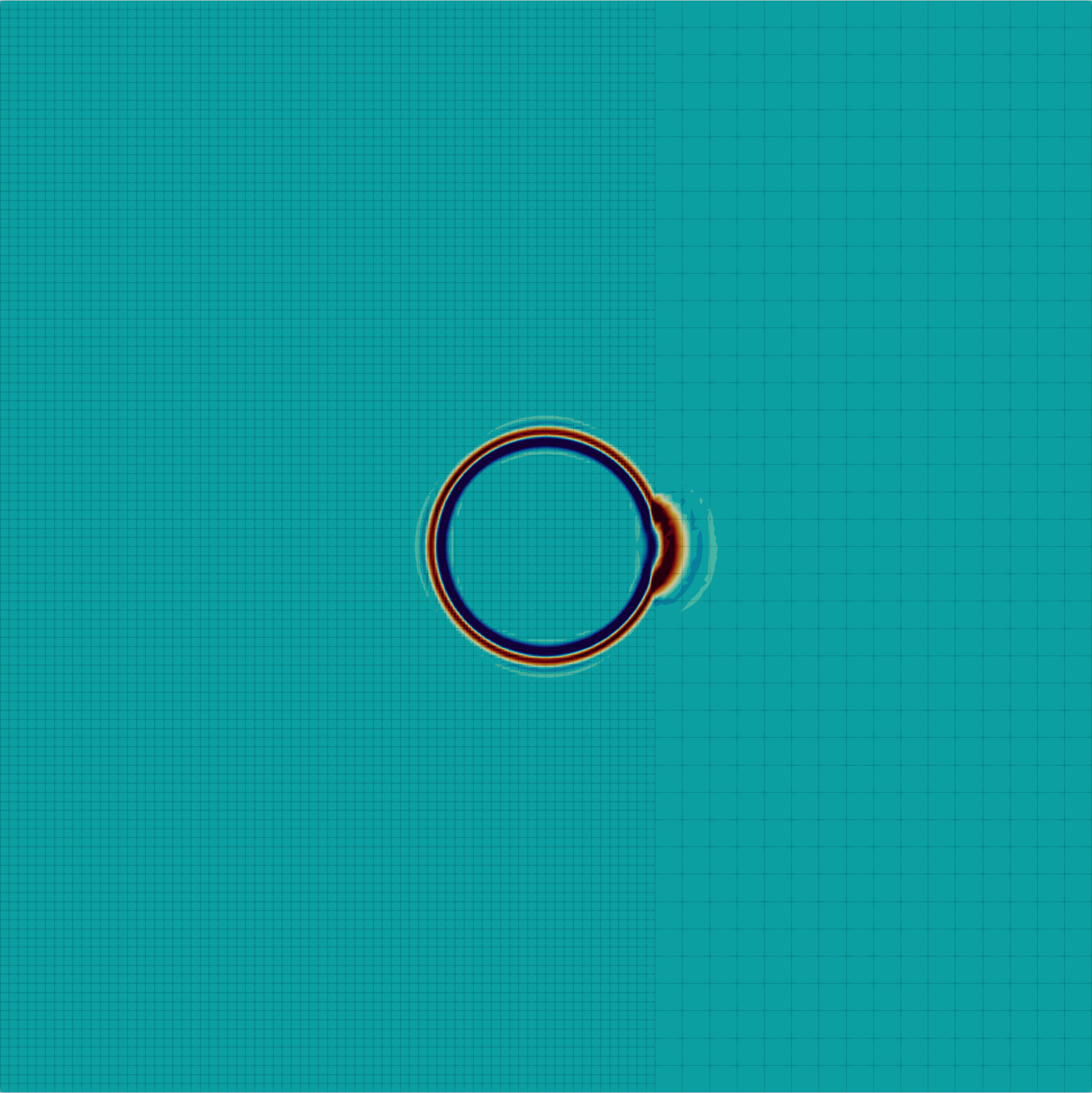}\label{fig:appl2}}
  
  \subfloat[][$x_1$-$x_2$ plane, $t=0.3\unit{s}$.]{\includegraphics[width=0.48\textwidth]{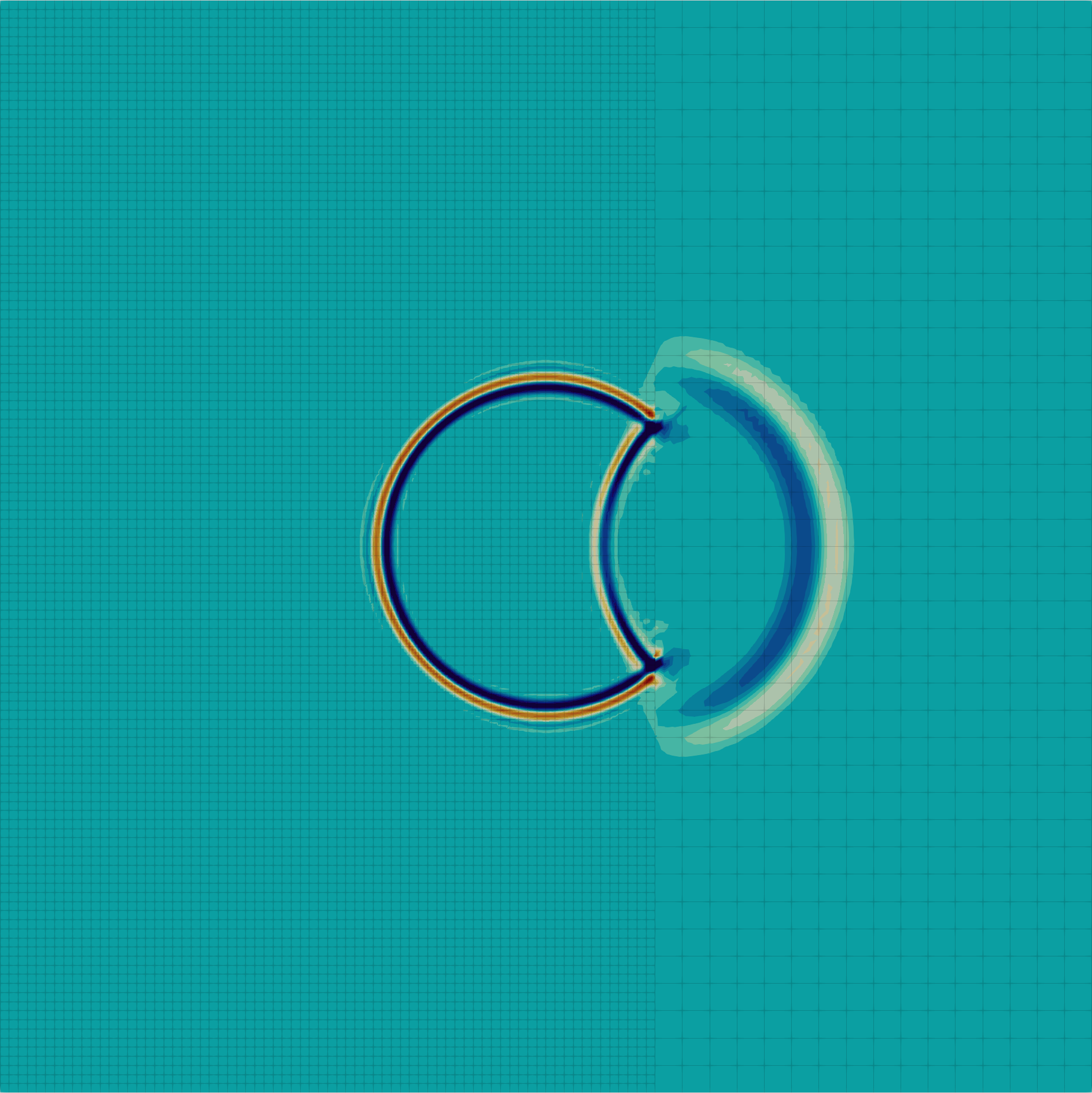}\label{fig:appl3}}\hfill
  \subfloat[][3D view, $t=0.38\unit{s}$.]{\includegraphics[width=0.48\textwidth]{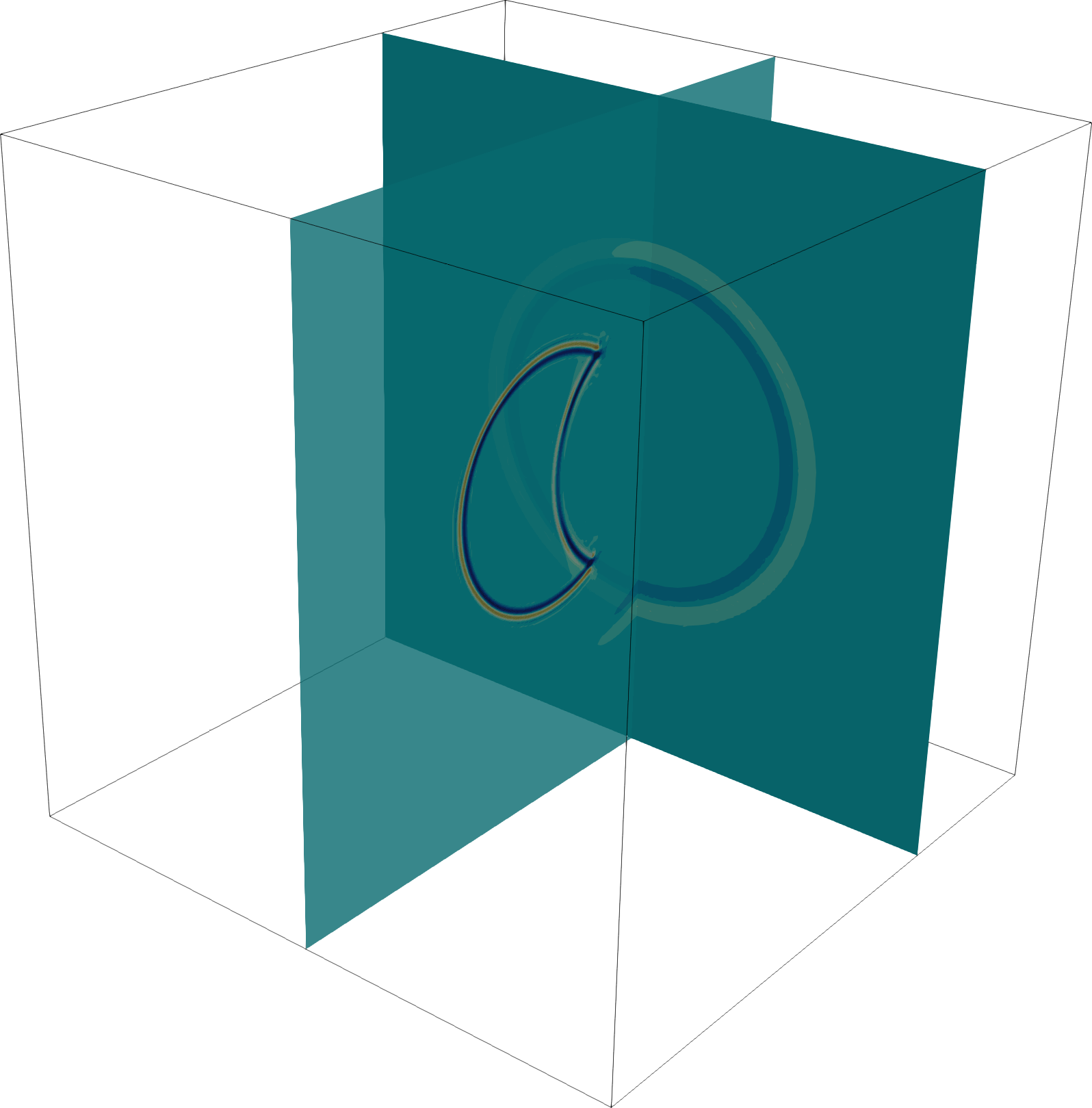}\label{fig:appl4}}\hfill
 \caption{Application: Snapshot of acoustic pressure at different times. At $t=0.3\unit{s}$ the transmitted, the reflected and the Huygens wave can be seen.\label{fig:appl}}
\end{figure*}

In the right domain, we use element sizes that are three times as big compared to the left domain to resolve both domains up to the same frequency.
In the left domain we use elements with maximum edge length $h_{\max,\Omega_1}=0.0167$ and accordingly $h_{\max,\Omega_2}=0.05$.
The used polynomial degree is $k=3$. 
The pressure field at different times can be seen in Figure~\ref{fig:appl}.

To quantify the effect of the NCI we also run the simulations on a domain with $h_{\max,\Omega_1}=h_{\max,\Omega_2}=0.0167$ and $h_{\max,\Omega_1}=h_{\max,\Omega_2}=0.05$.
We record the pressure at $1000$ points along $x_1$, $x_2=0$, $x_3=0$ at $t=0.2\unit{s}$.
The discretization with the smallest mesh size $h_{\max,\Omega_1}=h_{\max,\Omega_2}=0.0167$ serves as a reference.
It is supposed to produce the most accurate solution but uses too many DoFs if we want to resolve the same frequencies in both fluids.

\begin{figure*}[t]
  \centering
  \subfloat[][Pressure along the whole domain.]{\includegraphics[width=0.5\textwidth]{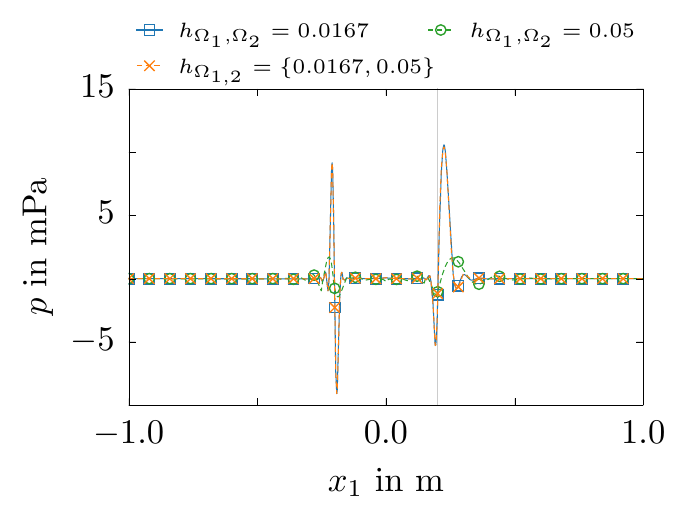}\label{fig:press}}\hfill
  \subfloat[][Detail view around NCI.]{\includegraphics[width=0.5\textwidth]{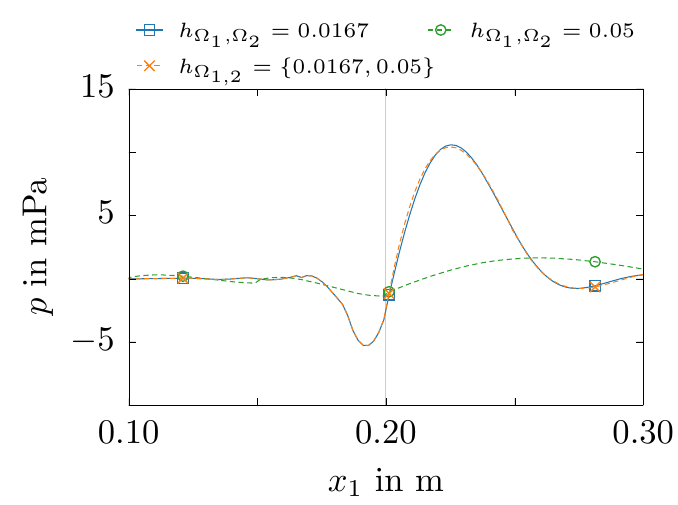}\label{fig:presszoom}}
\caption{Application: Pressure values along $x_1$ at $x_2=0, x_3=0$ and at $t=0.2\unit{s}$. The position of the NCI is indicated by the vertical line.\label{fig:lineplot}}  
\end{figure*}

The recorded pressure profile is plotted in Figure~\ref{fig:press}.
Figure~\ref{fig:presszoom} shows a detailed view around the interface.
We observe great differences to the reference for the discretization with the biggest mesh size $h_{\max,\Omega_1}=h_{\max,\Omega_2}=0.05$.
However, using the biggest and smallest mesh size for the different regions, employing the non-conforming formulation, gives a result that is in good agreement with the reference solution.
In this case, the finest domain has $442,368,000$ DoFs while the domain with different element sizes has $271,974,400$ DoFs, i.e. the problem size is reduced approximately by 40\% in comparison the fine problem while keeping the same accuracy.
This highly encourages to use NCIs for this kind of problems to effectively reduce the number of DoFs.
 
\section{Conclusion}
Using Nitsche-type mortaring, we proposed a stable non-conforming DG discretization for the acoustic conservation laws.
We showed that point-to-point interpolation is unsuitable in this setting since it introduces errors related to non-smooth representations of values in quadrature rules.
Therefore, we can not avoid the expensive computations of element intersections between primary and secondary elements.

The proposed method collects integration rules on the intersections between secondary volume elements and facets of primary elements.
This way, the method naturally extends to overlapping elements and is a perfect starting point for problems with rotating interfaces.
The method is subject to optimal spatial convergence rates.
Measuring the error region-wise, we can show that the method converges optimally in all sub-domains.
Nevertheless, errors are propagating in the domain;
therefore, optimal spatial convergence can only be applied in a meaningful way if triangulations are constructed such that errors are of the same magnitude in all parts of the domain.
Thus, we recommend using element sizes that resolve the same frequencies in all sub-domains in acoustics.
With an application, we demonstrated that this procedure efficiently reduces needed degrees of freedom while maintaining accuracy.

\subsection*{Acknowledgements}  
This project has received funding from the European Union's Framework Programme for Research and Innovation Horizon 2020 (2014-2020) under the Marie Sk\l{}odowska--Curie Grant Agreement No. [812719].\\
The computational results presented have been achieved in part using the Vienna Scientific Cluster (VSC).\\
The authors acknowledge collaboration with Marco Feder, Niklas Fehn, Martin Kronbichler and Magdalena Schreter, as well as the \texttt{deal.II} community.

\end{document}